\documentclass[11pt]{article}
\usepackage{amsmath, amsthm,amsfonts,amssymb,amscd}
\usepackage{color}
\usepackage{geometry}
\usepackage{tikz}
\usepackage{subfigure}
\numberwithin{equation}{section}
\usepackage[T1]{fontenc}
\usepackage{cite}
\usepackage{enumerate}

\usepackage{float}
\usepackage{soul}
\usepackage{marginnote}
\def\tto{\;{\lower 1pt \hbox{$\rightarrow$}}\kern -10pt
\hbox{\raise 2pt \hbox{$\rightarrow$}}\;}

\def\ra{\rangle}
\def\la{\langle}

\def\B{\mathbb B}

\def\h{\hfill\Box}
\def\R{\mathbb R}

\def\ox{\bar{x}}

\def\oz{\bar{z}}
\def\ov{\bar{v}}
\def\ou{\bar{u}}

\def\h{\hfill\triangle}

\def\ph{\varphi}

\newcounter{lk}

\begin{document}

\newtheorem{Theorem}{Theorem}[section]
\newtheorem{Proposition}[Theorem]{Proposition}
\newtheorem{Remark}[Theorem]{Remark}
\newtheorem{Lemma}[Theorem]{Lemma}
\newtheorem{Corollary}[Theorem]{Corollary}
\newtheorem{Definition}[Theorem]{Definition}
\newtheorem{Example}[Theorem]{Example}
\newtheorem{Fact}[Theorem]{Fact}
\renewcommand{\theequation}{\thesection.\arabic{equation}}
\normalsize
\def\proof{
\normalfont
\medskip
{\noindent\itshape Proof.\hspace*{6pt}\ignorespaces}}
\def\endproof{$\h$ \vspace*{0.1in}}

\title{\small \bf SOME RESULTS ON THE STRICT FRÉCHET DIFFERENTIABILITY OF THE METRIC PROJECTION OPERATOR IN HILBERT SPACES}

\date{}
\author{ LE VAN HIEN \footnote{Faculty of Pedagogy, Ha Tinh University, Ha Tinh city, Ha Tinh, Vietnam; email: hien.levan@htu.edu.vn}}
\maketitle

{\small \begin{abstract}
In this paper, we first present a simpler proof of a result on the strict Fréchet differentiability of the metric projection operator onto closed balls centered at the origin in Hilbert spaces, which given by  Li in \cite{Li24}. Then, based on this result, we prove the strict Fréchet differentiability of the metric projection operator onto closed balls with center at arbitrarily given point  in Hilbert spaces. Finally, we study the strict Fréchet differentiability of the metric projection operator onto the second-order cones in Euclidean spaces.

\end{abstract}}

\medskip

\normalsize
\section{Introduction}
\setcounter{equation}{0}

For a nonempty closed subset $C$ in a Banach space $X$, the distance function $d_C$
associated with $C$ is defined by
$$d_C(x)=\inf\{\|x-c\|:c\in C\},\ \ \mbox{for all}\ \ x\in X,$$
and the metric projection $P_C$ is given by
$$P_C(x)=\inf\{\ox\in C:d_C(x)=\|\ox-c\|\},\ \ \mbox{for all}\ \ x\in X.$$

The metric projection $P_C(x)$ is a set-valued mapping in the general case.  However, if $X$ is a Hilbert space and $C$ a nonempty closed convex subset, the metric projection operator is a well-defined single-valued mapping. Moreover, the metric projection operator enjoys properties that make it a useful and powerful tool in the optimization, computational mathematics, theory of equation, control theory, and others. Therefore, studying the properties of metric projection operators in Hilbert space is an interesting topic. The problem has been studied by many researchers for a long time, and many interesting results in this direction have been established; see \cite{FP82, H77, K84, Li23, Li24, LLX23, M02, N95, S87} and the references therein.

In this paper, we are interested in the strict Fréchet differentiability of the metric projection operator $P_C$, where $C$ is the closed ball in Hilbert spaces or a second-order cone in Euclidean spaces. For $C$  is a closed ball in real uniformly convex and uniformly smooth Banach space, recently, Li \cite{Li23}  has given the exact representations of the Gâteaux directional derivatives of $P_C$. In Hilbert space, the Gâteaux directional derivative of the metric projection operator onto the closed balls and their properties were investigated by Li et al. in \cite{LLX23}. Most recently, in \cite{Li24}, Li has investigated the strict Fréchet differentiability of the metric projection operator onto closed balls $r\B$ centered at the origin in Hilbert spaces. In his work, Li has proven strict Fréchet differentiability of $P_{r\B}$ and also found the exact expressions for Fréchet derivatives. In addition, he proved on $\mbox{bd}r\B$ that the Gâteaux directional differentiability of $P_{r\B}$ at $\ox$ does not imply its Fréchet differentiability at this point. However, for the metric projection operator onto closed balls with center at arbitrarily given point $c$ in the spaces, whether the results in \cite{Li24} are correct or not still remains unknown.

For $C=\mathcal{K}$, which is second-order cone in Euclidean spaces, Koranyi \cite{K84} has shown that $P_\mathcal{K}$ is Fréchet differentiable at every point not on the boundary of $\mathcal{K}$ and $-\mathcal{K}$, and gave the Fréchet derivative formula of the metric projection operator in this case. Using the results of Fréchet differentiability of $P_\mathcal{K}$ in \cite{K84}, Outrata and  Sun~\cite{OS08} established a formula for exactly computing the the limiting (Mordukhovich) coderivative of $P_\mathcal{K}$. Also, they gave some applications of the solution map of the parameterized second-order cone complementarity problem and the mathematical program with a second-order cone complementarity problem among the constraints.

The aim of this paper is to study the strict Fréchet differentiability of the metric projection operator onto closed balls centered at arbitrarily given point  in Hilbert spaces. More precisely, we first present a simpler proof of a result about the strict Fréchet differentiability of the metric projection operator $P_{r\B}$ with $r\B=\B(0,r)=\{x\in X| \|x\|\leq r\}$, which given in \cite{Li24}. Then, based on this result, we prove the strict Fréchet differentiability of the metric projection operator $P_{\B(c,r)}$ with $\B(c,r)=\{x\in X|\|x-c\|\leq r\}$ for any $c\in X.$ Finally, we study the strict Fréchet differentiability of the metric projection operator $P_\mathcal{K}$, where $\mathcal{K}=\{z=(z_1,z_2)\in\R\times\R^{m-1}||z_1|>\|z_2\|\}$ is a second-order cone in $\R^m$.

The structure of the paper is as follows.   In section 2,  we recall some preliminary materials. In section 3, we present results on strict Fréchet differentiability of the metric projection operator onto closed balls in Hilbert spaces. Then, in section 4, we study the strict Fréchet differentiability of the metric projection operator onto the second-order cones in Euclidean spaces. Finally, we conclude the paper in section 5 where we discuss some perspectives of the obtained results and future works.

\section{Preliminaries}
\setcounter{equation}{0}
 
  We first give the following notations that will be used throughout the paper. Let $X$ be Banach spaces, we denote by $\|.\|$ the norm in $X$. Given a set $C \subset X,$ we denote by $C^o$, $\mbox{bd}C$ its interior, boundary respectively. The closed ball with center $\ox$ and radius $r>0$ is denoted by $\B(\ox,r)$.  For a function $g:X\to\R$, we denote $\big(g(x)\big)_+=\max(0, g(x))$ and $\big(g(x)\big)_-=\min(0, g(x))$. Define the identity mapping on $X$ by $I_X:X\to X$ with $I_X(x)=x,\ \mbox{for all}\ x\in X.$  Denote $t\downarrow 0$ means that $t\to 0$ and  $t>0.$


\begin{Definition}{\rm
Let $X, Y$ be Banach spaces,  $\Omega\subset X$. The single-valued mapping $f:X\to Y$ is called {\it Lipschitz continuous} on $\Omega$ if there exists a real constant $L\geq 0$ such that
$$\|f(x_1)-f(x_2)\|\leq L\|x_1-x_2\|,\ \ \mbox{ for all}\ x_1, x_2 \in \Omega.$$}
\end{Definition}
Next we recall some concepts about "the differentiability" of the single-valued mappings in Banach spaces.  
\begin{Definition}{\rm(see \cite{Li24})
Let $X, Y$ be Banach spaces, $f:X\to Y$ be a single-valued mapping. For $x, w\in X$ with $w\not=0$, if the
following limit exists, which is a point in $Y$,
$$\lim\limits_{t\downarrow0}\dfrac{f(x+tw)-f(x)}{t}:=f'(x)(w),$$
then, $f$ is said to be {\it Gâteaux directionally differentiable} at point $x$ along direction $w$. $f'(x)(w)$ is
called the {\it Gâteaux directional derivative} of $f$ at point $x$ along direction $w$. 

Let $A\subset X,$ if $f$ is Gâteaux directionally differentiable at every point $x\in A$ then $f$ is said to be Gâteaux directionally differentiable on $A$.

}
\end{Definition}
\begin{Definition}
	{\rm (see \cite{Li24})
	 Let $X, Y$ be Banach spaces, $f:X\to Y$ be a single-valued mapping. For any given $\ox\in X$ if there is a linear continuous mapping $\nabla f(\ox):X\to Y$ such that
		$$\lim\limits_{h\to 0}\dfrac{\|f(\ox+h)-f(\ox)-\nabla f(\ox)(h)\|}{\|h\|}=0$$
	then $f$ is said to be {\it Fréchet differentiable} at $\ox$ and $\nabla f(\ox)$ is called the {\it Fréchet
derivative} of $f$ at $\ox$.

More strictly, if
	$$\lim\limits_{(u,v)\to (\ox,\ox)}\dfrac{\|f(u)-f(v)-\nabla f(\ox)(u-v)\|}{\|u-v\|}=0$$
then $f$ is said to be {\it  strictly Fréchet differentiable} at $\ox$.

We say $f$ is {\it continuously differentiable} at $\ox$ if its differentiability on a neighborhood $U$ of $\ox$ and $\nabla f(.)$ is continuous at $\ox$.

Let $A\subset X,$ if $f$ is Fréchet differentiable (strictly Fréchet differentiable, continuously differentiable) at every point $x\in A$ then $f$ is said to be Fréchet differentiable (strictly Fréchet differentiable, continuously differentiable, respectively) on $A$.}\end{Definition}

\begin{Remark} {\rm
(i) $f$ is strictly Fréchet differentiable at $\ox$

\hspace{2.1cm}$\Longrightarrow$ $f$ is  Fréchet differentiable at $\ox$

\hspace{2cm} $\Longrightarrow$ $f$ Gâteaux directionally differentiable at $\ox$ along every direction, and
$$f'(\ox)(w)=\nabla f(\ox)(w),\ \ \mbox{for all} \ \ w\in X\backslash\{0\}.$$
(ii) $f$ Gâteaux directionally differentiable at $\ox$ along every direction

$\not\Longrightarrow$ $f$ is  Fréchet differentiable at $\ox$. (see Proposition \ref{Pro32} and Theorems \ref{Thm34}, \ref{Thm42}.)
}\end{Remark}

\begin{Proposition}\label{Pro23}
	Let $X, Y$ be Banach spaces, $f,g: X\to Y$ are Fréchet differentiable (strictly Fréchet differentiable, continuously differentiable) at $\ox\in X.$ Then, $\alpha f+\beta g$ is Fréchet differentiable (strictly Fréchet differentiable, continuously differentiable, respectively) at $\ox$, and 
	$$\nabla (\alpha f+\beta g)(\ox)=\alpha\nabla f(\ox)+\beta\nabla g(\ox),\ \ \mbox{for all} \ \alpha, \beta\in \R.$$
\end{Proposition}

\begin{Proposition}\label{Pro24} {\rm (see \cite[Theorem 2.1]{A99})}
	Let $X, Y, Z$ be Banach spaces, $f: X\to Y$ is Fréchet differentiable (strictly Fréchet differentiable, continuously differentiable) at $\ox\in X,$ $g:Y\to Z$ is Fréchet differentiable (strictly Fréchet differentiable, continuously differentiable) at $f(\ox).$ Then, $g\circ f$ is Fréchet differentiable (strictly Fréchet differentiable, continuously differentiable, respectively) at $\ox$, and 
	$$\nabla (g\circ f)(\ox)=\nabla g(f(\ox))\circ\nabla f(\ox).$$
\end{Proposition}

{\bf Proof.}
Suppose $f: X\to Y$ is strictly Fréchet differentiable,  at $\ox\in X,$ $g:Y\to Z$ is stritcly Fréchet differentiable at $f(\ox).$

Since $g$ is strictly Fréchet differentiable at $f(\ox)$, for any $\epsilon>0,$ there exists $r>0$ such that
$$\|g(s)-g(t)-\nabla g(f(\ox))(s-t)\|\leq \epsilon\|s-t\|, \ \ \mbox{for all } \ s,t\in \B(f(\ox), r).$$
Since $f$ is strictly Fréchet differentiable at $\ox$, there exists $\delta$ such that
$$\|f(u)-f(v)-\nabla f(\ox)(u-v)\|\leq \epsilon\|u-v\| \ \mbox{and}\ \ f(u),f(v)\in\B(f(\ox),r), \ \ \mbox{for all } \ u,v\in \B(\ox, \delta).$$
In what follows we assume that $u,v\in \B(\ox, \delta).$ Then we have
$$\|f(u)-f(v)\|\leq \|f(u)-f(v)-\nabla f(\ox)(u-v)\|+\|\nabla f(\ox)(u-v)\|\leq \big(\epsilon+\|\nabla f(\ox)\|\big)\|u-v\|.$$
Moreover,
$$\begin{array} {rl}&\|\nabla g(f(\ox))(f(u)-f(v))-\nabla g(f(\ox))(\nabla f(\ox)(u-v))\|\\&\leq \|\nabla g(f(\ox))\|\|f(u)-f(v)-\nabla f(\ox)(u-v)\|\\&\leq \epsilon\|\nabla g(f(\ox))\|\|u-v\|.\end{array}$$
Since $f(u),f(v)\in\B(f(\ox),r)$, we have 
$$\|g(f(u))-g(f(v))-\nabla g(f(\ox))(f(u)-f(v))\|\leq \epsilon\|f(u)-f(v)\|\leq \epsilon\big(\epsilon+\|\nabla f(\ox)\|\big)\|u-v\|.$$
Therefore
$$\begin{array}{rl} &
	\|g\circ f(u)-g\circ f(v)-\nabla g(f(\ox))\circ\nabla f(\ox)(u-v)\|\\&\leq \|g(f(u))-g(f(v))-\nabla g(f(\ox))(f(u)-f(v))\| \\&+ \|\nabla g(f(\ox))(f(u)-f(v))-\nabla g(f(\ox))(\nabla f(\ox)(u-v))\|\\
	&\leq \epsilon\big(\epsilon+\|\nabla f(\ox)\|\big)\|u-v\|+\epsilon\|\nabla g(f(\ox))\|\|u-v\|\\
	&=\epsilon\big(\epsilon+\|\nabla f(\ox)\|+\|\nabla g(f(\ox))\|\big)\|u-v\|.
	\end{array}$$
This shows that $g\circ f$ is strictly Fréchet differentiable at $\ox$, and 
$$\nabla (g\circ f)(\ox)=\nabla g(f(\ox))\circ\nabla f(\ox).$$
\hfill $\square$

 \section{Strict Fréchet differentiability of the metric projection operator onto balls}\label{Sec3}
\setcounter{equation}{0}

In this section, let $(H,\| . \|)$ be a real Hilbert space with inner product $\la\ .\ \ra$. For any $c\in H$ and $r>0$, let $\B(c,r)$ denote the closed ball in $H$ with  centered at $c$ and  radius $r$. Notation $\B(0,1)=\B$ and $\B(0,r)=r\B.$
Let $P_{\B(c,r)}:H\to \B(c,r)$ be a  metric projection operator onto  $\B(c,r)$ in $H$. It is well-known that $P_{\B(c,r)}$ is a well-defined single-valued mapping, and
$$P_{\B(c,r)}(x)=\begin{cases}
	x\ \ \ \ \ \ \ \ \ \ \ \ \ \ \ \ \ \ \ \ \ \ \ \ \  \mbox{if} \ x\in \B(c,r);\\
	c+\dfrac{r}{\|x-c\|}(x-c)\ \ \ \mbox{if} \ x\in H\backslash\B(c,r).
\end{cases}$$
In particular, 
$$P_{\B}(x)=\begin{cases}
	x\ \ \ \ \ \ \ \ \ \ \ \ \  \mbox{if} \ x\in \B;\\
	\dfrac{x}{\|x\|}\ \ \ \ \ \ \ \  \ \mbox{if} \ x\in H\backslash\B.
\end{cases}$$
For any $x\in\B(c,r)$, we define two subsets $x^\uparrow_{(c,r)}$ and $x^\downarrow_{(c,r)}$ of $H\backslash\{0\}$ as follows: 
$$x^\uparrow_{(c,r)}=\{v\in H\backslash\{0\}|\, \mbox{there is}\  \delta>0 \ \mbox{such that}\ \|(x+tv)-c\|\geq r,\ \mbox{for all}\ t\in (0,\delta)\}.$$
$$x^\downarrow_{(c,r)}=\{v\in H\backslash\{0\}|\, \mbox{there is}\  \delta>0 \ \mbox{such that}\ \|(x+tv)-c\|<r,\ \mbox{for all}\ t\in (0,\delta)\}.$$
In particular, 
$$x^\uparrow:=x^\uparrow_{0,1}=\{v\in H\backslash\{0\}|\, \mbox{there is}\  \delta>0 \ \mbox{such that}\ \|x+tv\|\geq r,\ \mbox{for all}\ t\in (0,\delta)\}.$$
$$x^\downarrow:=x^\downarrow_{(0,1)}=\{v\in H\backslash\{0\}|\, \mbox{there is}\  \delta>0 \ \mbox{such that}\ \|x+tv\|<r,\ \mbox{for all}\ t\in (0,\delta)\}.$$

For any $x\in H\backslash \{0\}$, let $S(x)$ be the one-dimensional subspace of $H$ generated by $x$. Let $O(x) $ denote the orthogonal subspace of $x$ in $H$. Then, we have
$$H=S(x)\oplus O(x).$$ For any fixed $x\in H\backslash\{0\}$, we define a real valued function $a(x, . ): H\to \R$ by
$$a(x, u)=\dfrac{\la x, u\ra}{\|x\|^2}, \ \ \ \ \mbox{for all}\ u\in H.$$
And a mapping  $o(x, . ):H\to O(x)$ by
$$o(x, u)=u-\dfrac{\la x, u\ra}{\|x\|^2}x, \ \ \ \ \mbox{for all}\ u\in H.$$
We see that \\

 $\blacklozenge$\ \ \ \ \  $u=a(x,u)x+o(x,u)\ \ \ \mbox{and}\ \ \ \big\la a(x,u)x, o(x,u)\big\ra=0, \ \ \mbox{for all}\ u\in H.$

$\blacklozenge$\ \ \ \ \ $a(x, . ), o(x, . )$ are continuous linear mappings, and $o(x,x)=0$, $o(x,x+h)=o(x,h)$ for all $h\in H.$ 

$\blacklozenge$\ \ \ \ \  $ \|o(x,h)\|\leq \|h\|+\dfrac{|\la x,h\ra|}{\|x\|^2}\|x\|\leq \|h\|\Big( 1+\dfrac{1}{\|x\|}\Big), \ \mbox{for all}\ h\in H.$

$\blacklozenge$\ \ \ \ \ $a(x,u)\to 1$ and $o(x,u)\to 0$ as $u\to x.$\\

To proceed, we need the following result, which shows that strict Fréchet  differentiability can be obtained from continuous differentiability.

\begin{Lemma}{\rm (see \cite[Theorem 25.23]{S96}) } \label{Pro31} Let $X$ and $Y$ be Banach spaces, let $\Omega\subset X$ be an open convex set, and let $f:X \to Y$ be Fréchet differentiable mapping on $\Omega$ with  $\nabla f(.)$ Lipschitz continuous on it. Then $f$ is strict Fréchet differentiable at $\ox$, for every $\ox\in\Omega$.
\end{Lemma} 
{\bf Proof.}
Suppose that $f$ is Fréchet differentiable and $\nabla f(.)$ Lipschitz continuous on $\Omega$ with modulus $L$.\\
Given $u,v\in \Omega$, define $\ph:[0,1]\to X$ by
$$\ph(t)=f(tu+(1-t)v).$$
Thanks to the continuous differentiability of $f$ on $\Omega$, by Proposition \ref{Pro24}, we obtain the continuous differentiability of $\ph$ on $[0,1]$, and
$$\nabla\ph(t)=\nabla f(tu+(1-t)v)(u-v).$$
Therefore
$$\ph(1)-\ph(0)=\int_{0}^{1}\nabla\ph(t)dt.$$
It means
$$f(u)-f(v)=\int_{0}^{1}\nabla f(tu+(1-t)v)(u-v)dt.$$
This implies
$$f(u)-f(v)-\nabla f(\ox)(u-v)=\int_{0}^{1}\big(\nabla f(tu+(1-t)v)-\nabla f(\ox)\big)(u-v)dt.$$
So
$$\begin{array}{rl}\|f(u)-f(v)-\nabla f(\ox)(u-v)\|&\leq \int_{0}^{1}\big\|\big(\nabla f(tu+(1-t)v)-\nabla f(\ox)\big)(u-v)\big\|dt\\
	&\leq \|u-v\|\int_{0}^{1}\big\|\nabla f(tu+(1-t)v)-\nabla f(\ox)\big\|dt\\
	&\leq \|u-v\|\int_{0}^{1}L\big\|(tu+(1-t)v)-\ox\big\|dt\\
		&\leq L\|u-v\|\int_{0}^{1}\big(t\|u-\ox\|+(1-t)\|v-x\|\big)dt\\
	&=\dfrac{L}{2} \big(\|u-\ox\|+\|v-x\|\big)\|u-v\|.
		\end{array}$$
Therefore, we have
$$\begin{array}{rl}0&\leq\lim\limits_{(u,v)\to(\ox,\ox)}\dfrac{\|f(u)-f(v)-\nabla f(\ox)(u-v)\|}{\|u-v\|}\\
	&\leq \lim\limits_{(u,v)\to(\ox,\ox)}\dfrac{L}{2} \big(\|u-\ox\|+\|v-x\|\big)\\&=0.\end{array}$$
This shows that  $f$ is strict Fréchet differentiable at $\ox$.\hfill$\square$\\

The following result has been given in \cite{Li24}, which provides information about the strict Fréchet differentiability of the metric projection onto a closed ball with centered at the origin. We present another proof, simpler than the proof in \cite{Li24}.
\begin{Proposition}{\rm (see \cite[Theorem 3.1]{Li24}) } \label{Pro32} 
Let $H$ be a Hilbert space. For any  $r>0$, the metric projection $P_{r\B}:H\rightarrow r\B$   has the following differentiability properties.

(i)  $P_{r\B}$ is strictly Fréchet differentiable on $r\B^o$ satisfying $$\nabla P_{r\B}(\ox)=I_H \ \ \mbox{for every}\ \ox\in r\B^o$$

(ii) $P_{r\B}$ is strictly Fréchet differentiable on $H\backslash r\B$ such that, for every $\ox\in H\backslash r\B$, 
$$\nabla P_{r\B}(\ox)(u)=\dfrac{r}{\|\ox\|}\Big(u-\dfrac{\la \ox, u\ra}{\|\ox\|^2}\ox\Big)\ \ \mbox{for all} \ u\in H.$$

(iii) For the subset $\mbox{bd}r\B$, we have

(a)  $P_{r\B}$  is Gâteaux directional differentiable on $\mbox{bd}r\B$ satisfying that, for every point  $\ox\in \mbox{bd}r\B$, the following representations are satisfied\\
$$P^{'}_{r\B}(\ox)(w)=\begin{cases}
	w-\dfrac{1}{r^2}\la \ox,w\ra\ox \ \ \ \mbox{if} \ \ \ w\in\ox^\uparrow_r;\\
	0 \ \ \ \ \ \ \ \ \ \ \ \ \ \ \ \ \ \ \ \mbox{if}\ \  \ w=\ox;\\w \ \ \ \  \ \ \ \ \ \ \ \ \ \  \ \ \ \  \mbox{if} \ \ \ \ w\in\ox^\downarrow_r.
\end{cases}$$  

(b)  $P_{r\B}$ is not Fréchet differentiable at any point $\ox\in \mbox{bd}r\B$. That is, $\nabla P_{r\B}(\ox)$ 
does not exist, for any $\ox\in \mbox{bd}r\B$. 
\end{Proposition}

{\bf Proof.}
Proof of (ii). Since $\ox\in H\backslash r\B$, exits $\delta>0$ such that $\B(\ox,\delta)\subset H\backslash r\B$. For any $x\in \B(\ox,\delta)$, we prove that $P_{r\B}$ is Fréchet differentiable at $x$, and
$$\nabla P_{r\B}(x)(h)=\dfrac{r}{\|x\|}\Big(h-\dfrac{\la x, h\ra}{\|x\|^2}x\Big)=\dfrac{r}{\|x\|}o(x,h)\ \ \mbox{for all} \ h\in H. $$
We have 
$$\begin{array} {rl}&\lim\limits_{h\to 0}\dfrac{P_{r\B}(x+h)-P_{r\B}(x)-\dfrac{r}{\|x\|}o(x,h)}{\|h\|}
	\\&=\lim\limits_{h\to 0}\dfrac{\dfrac{r(x+h)}{\|x+h\|}-\dfrac{rx}{\|x\|}-\dfrac{r}{\|x\|}o(x,h)}{\|h\|}\\
&=r\lim\limits_{h\to 0}\dfrac{\dfrac{a(x,x+h)+o(x,x+h)}{\|a(x,x+h)+o(x,x+h)\|}-\dfrac{x}{\|x\|}-\dfrac{o(x,h)}{\|x\|}}{\|h\|}\\
&=r\lim\limits_{h\to 0}\dfrac{\dfrac{a(x,x+h)}{\|a(x,x+h)+o(x,x+h)x\|}-\dfrac{x}{\|x\|}}{\|h\|}+r\lim\limits_{h\to 0}\dfrac{\dfrac{o(x,x+h)}{\|x+h\|}-\dfrac{o(x,h)}{\|x\|}}{\|h\|}\\\\
\end{array}$$
So,
$$\begin{array} {rl}&
	\bigg\|\dfrac{a(x,x+h)x}{\|a(x,x+h)x+o(x,x+h)\|}-\dfrac{x}{\|x\|}\bigg\|\\&=\dfrac{\big|a(x,x+h)\|x\|-\|a(x,x+h)x+o(x,x+h)\|\big|}{\|a(x,x+h)x+o(x,x+h)\|}\\
	&=\dfrac{\big|(a(x,x+h))^2\|x\|^2-\|a(x,x+h)x+o(x,x+h)\|^2\big|}{\|a(x,x+h)x+o(x,x+h)\|\big|a(x,x+h)\|x\|+\|a(x,x+h)x+o(x,x+h)\|\big|}\\
	&=\dfrac{\|o(x,x+h)\|^2}{\|a(x,x+h)x+o(x,x+h)\|\big|a(x,x+h)\|x\|+\|a(x,x+h)x+o(x,x+h)\|\big|}\\
	&\leq \dfrac{\|h\|^2\Big( 1+\dfrac{1}{\|x\|}\Big)^2}{\|a(x,x+h)x+o(x,x+h)\|\big|a(x,x+h)\|x\|+\|a(x,x+h)x+o(x,x+h)\|\big|}
\end{array}$$
This means,
$$ \begin{array} {rl}

0&\leq \dfrac{\bigg\|\dfrac{a(x,x+h)}{\|a(x,x+h)+o(x,x+h)x\|}-\dfrac{x}{\|x\|}\bigg\|}{\|h\|}\\
&\leq \dfrac{\|h\|\Big( 1+\dfrac{1}{\|x\|}\Big)^2}{\|a(x,x+h)x+o(x,x+h)\|\big|a(x,x+h)\|x\|+\|a(x,x+h)x+o(x,x+h)\|\big|}\\
&\to 0, \ \ \mbox{as}\ \ h\to 0.\end{array}$$
Therefore,
\begin{equation}\label{31}\lim\limits_{h\to 0}\dfrac{\dfrac{a(x,x+h)}{\|a(x,x+h)+o(x,x+h)x\|}-\dfrac{x}{\|x\|}}{\|h\|}=0.
\end{equation}
Moreover, we have
$$\begin{array} {rl}0&\leq \bigg\|\dfrac{o(x,x+h)}{\|x+h\|}-\dfrac{o(x,h)}{\|x\|}\bigg\|\\
	&=\|o(x,h)\|\bigg|\dfrac{1}{\|x+h\|}-\dfrac{1}{\|x\|}\bigg|\\
	&=\dfrac{\|o(x,h)\|\big|\|x\|-\|x+h\|\big|}{\|x\|\|x+h\|}\\
&\leq \dfrac{\|h\|\bigg( 1+\dfrac{1}{\|x\|}\bigg)\|x - (x+h)\|}{\|x\|\|x+h\|}\\
&=\dfrac{\|h\|^2\bigg( 1+\dfrac{1}{\|x\|}\bigg)}{\|x\|\|x+h\|}.
		\end{array}$$
This means, $$0\leq \dfrac{\bigg\|\dfrac{o(x,x+h)}{\|x+h\|}-\dfrac{o(x,h)}{\|x\|}\bigg\|}{\|h\|}\leq \dfrac{\|h\|\bigg( 1+\dfrac{1}{\|x\|}\bigg)}{\|x\|\|x+h\|}\to 0 \ \ \mbox{as}\ h\to 0.$$
So, 
 \begin{equation}\label{32}\lim\limits_{h\to 0}\dfrac{\dfrac{o(x,x+h)}{\|x+h\|}-\dfrac{o(x,h)}{\|x\|}}{\|h\|}=0.\end{equation}
By \eqref{31}-\eqref{32}, we get 
$$\lim\limits_{h\to 0}\dfrac{P_{r\B}(x+h)-P_{r\B}(x)-\dfrac{r}{\|x\|}o(x,h)}{\|h\|}=0.$$
This proves that $P_{r\B}$ is Fréchet differentiable at $x$, and
$$\nabla P_{r\B}(x)(h)=\dfrac{r}{\|x\|}o(x,h)\ \ \mbox{for all} \ h\in H. $$

Next we prove $\nabla P_{r\B}(.)$ is Lipschitz continuous function  on $B(\ox, \delta)$. Indeed, for any $x_1, x_2\in\B(\ox, \delta)$, $u\in H$, since
$$\begin{array} {rl}&\Big(\nabla P_{r\B}(x_1)-\nabla P_{r\B}(x_2)\Big)(u)\\
&	=\dfrac{r}{\|x_1\|}\bigg(u-\dfrac{\la x_1,u\ra}{\|x_1\|^2}x_1\bigg)- \dfrac{r}{\|x_2\|}\bigg(u-\dfrac{\la x_2,u\ra}{\|x_2\|^2}x_2\bigg)\\
&=r\bigg(\dfrac{1}{\|x_1\|} -\dfrac{1}{\|x_2\|}\bigg)u-r\bigg(\dfrac{\la x_1,u\ra}{\|x_1\|^3}x_1-\dfrac{\la x_2,u\ra}{\|x_2\|^3}x_2  \bigg)\\
&=r\bigg(\dfrac{1}{\|x_1\|} -\dfrac{1}{\|x_2\|}\bigg)u-r\bigg(\dfrac{\la x_1,u\ra}{\|x_1\|^3}x_1-\dfrac{\la x_1,u\ra}{\|x_2\|^3}x_2 \bigg)-r\bigg(\dfrac{\la x_1,u\ra}{\|x_2\|^3}x_2-\dfrac{\la x_2,u\ra}{\|x_2\|^3}x_2  \bigg)\\
&=\dfrac{r(\|x_2\|-\|x_1\|)}{\|x_1\|\|x_2\|}u-\dfrac{r\la x_1,u\ra \big(\|x_2\|^3x_1-\|x_1\|^3x_2\big)}{\|x_1\|^3\|x_2\|^3}-\dfrac{r}{\|x_2\|^3}\la x_1-\ox,u\ra x_2.
\end{array}$$
We have 
$$\begin{array} {rl}&\Big\|\Big(\nabla P_{r\B}(x_1)-\nabla P_{r\B}(x_2)\Big)(u)\Big\|\\&=\bigg\|\dfrac{r(\|x_2\|-\|x_1\|)}{\|x_1\|\|x_2\|}u-\dfrac{r\la x_1,u\ra \big(\|x_2\|^3x_1-\|x_1\|^3x_2\big)}{\|x_1\|^3\|x_2\|^3}-\dfrac{r}{\|x_2\|^3}\la x_1-x_2,u\ra x_2\bigg\|\\
&	\leq \dfrac{r\big|\|x_2\|-\|x_1\|\big|}{\|x_1\|\|x_2\|}\|u\|+\dfrac{r|\la x_1,u\ra| \big\|\|x_2\|^3x_1-\|x_1\|^3x_2\big\|}{\|x_1\|^3\|x_2\|^3}+\dfrac{r}{\|x_2\|^3}|\la x_1-x_2,u\ra|\| x_2\|\\
&\leq \dfrac{r\big|\|x_2\|-\|x_1\|\big|}{\|x_1\|\|x_2\|}\|u\|
+\dfrac{r\|u\| \|x_1-x_2\|}{\|x_1\|^2}
+\dfrac{r\|u\| \big|\|x_2\|^3-\|x_1\|^3\big|}{\|x_1\|^2\|x_2\|^2}
+\dfrac{r}{\|x_2\|^2}\|x_1-x_2\|\|u\|\\
&=r\Bigg(\dfrac{\big|\|x_2\|-\|x_1\|\big|}{\|x_1\|\|x_2\|}+\dfrac{ \|x_1-x_2\|}{\|x_1\|^2} + \dfrac{ \big|\|x_2\|^3-\|x_1\|^3)\big|\|}{\|x_1\|^2\|x_2\|^2}
+\dfrac{\|x_1-x_2\|}{\|x_2\|^2}\Bigg)\|u\|.
\end{array}$$
So, 
$$\begin{array} {rl}0&\leq\Big\|\Big(\nabla P_{r\B}(x_1)-\nabla P_{r\B}(x_2)\Big)\Big\|\\
	&\leq r\Bigg(\dfrac{\big|\|x_2\|-\|x_1\|\big|}{\|x_1\|\|x_2\|}+\dfrac{ \|x_1-x_2\|}{\|x_1\|^2} + \dfrac{ \big|\|x_2\|^3-\|x_1\|^3)\big|\|}{\|x_1\|^2\|x_2\|^2}
	+\dfrac{\|x_1-x_2\|}{\|x_2\|^2}\Bigg)
\\
	&\leq r\Bigg(\dfrac{\|x_1-x_2\|}{\|x_1\|\|x_2\|}
	+\dfrac{ \|x_1-x_2\|}{\|x_1\|^2} 
	+ \dfrac{ \|x_1-x_2\|\big(\|x_2\|^2+\|x_2\|\|x_1\|+\|x_1\|^2\big)}{\|x_1\|^2\|x_2\|^2}
+\dfrac{\|x_1-x_2\|}{\|x_2\|^2}\Bigg)
\\
	&= 2r\Bigg(\dfrac{ 1}{\|x_1\|^2} 
+ \dfrac{ 1}{\|x_1\|\|x_2\|}
+\dfrac{1}{\|x_2\|^2}\Bigg)\|x_1-x_2\|\\
&\leq\dfrac{6r}{(\|\ox\|-\delta)^2}\|x_1-x_2\|.
\end{array}$$

This shows that $\nabla P_{r\B}(.)$ is Lipschitz continuous function  on $B(\ox, \delta)$.  By Lemma \ref{Pro31}, $P_{r\B}$ is strict Fréchet differentiable at $\ox$.

Proof of (iii.b). Suppose to contrary that $P_{r\B}$ is Fréchet differentiable at $\ox \in\mbox{bd}r\B.$ Then, there is a linear continuous mapping $\nabla P_{r\B}$, such that $$\nabla P_{r\B}(\ox)(w)=P^{'}_{r\B}(\ox)(w), \ \ \mbox{for all} \ w\in H.$$
By (iii.a), we have 
\begin{equation}\label{13}\nabla  P_{r\B}(\ox)(\ox) = P^{'}_{r\B}(\ox)(\ox)=0.\end{equation}
And $$\nabla  P_{r\B}(\ox)(-\ox) = P^{'}_{r\B}(\ox)(-\ox)=-\ox.$$
Thanks for the linearity of $\nabla  P_{r\B}(\ox)$, we get $\nabla  P_{r\B}(\ox)(\ox)=\ox$. This contradicts to \eqref{13}, which proves that $P_{r\B}$ is not Fréchet differentiable at $\ox.$
\hfill $\square$\\

\begin{Lemma} \label{Lem33}
	Let $C=\B(c,r)$, where $c\in H, r>0$. Then, we have
	$$P_C(x)=P_{r\B}(x-c)+c,\ \ \ \mbox{for all}\ x\in H.$$
\end{Lemma}

	

We now arrive at the main result of this section, the following theorem is to prove the strict Fréchet differentiability of the metric projection onto the closed balls with center at arbitrarily given point $c$ in $H$.
\begin{Theorem} \label{Thm34}
Let $H$ be a Hilbert space. For any  $c\in H, r>0$ put $C=\B(c,r)$. Then, the metric projection $P_C:H\rightarrow C$   has the following differentiability properties.

(i)  $P_C$ is strictly Fréchet differentiable on $C^o$ satisfying $$\nabla P_{C}(\ox)=I_H \ \ \mbox{for every}\ \ox\in C^o$$.

(ii) $P_C$ is strictly Fréchet differentiable on $H\backslash C$ such that, for every $\ox\in H\backslash C$, 
$$\nabla P_{C}(\ox)(u)=\dfrac{r}{\|\ox-c\|}\Big(u-\dfrac{\la \ox-c, u\ra}{\|\ox-c\|^2}(\ox-c)\Big)\ \ \mbox{for all} \ u\in H.$$

(iii) For the subset $\mbox{bd}C$, we have

(a)  $P_{C}$  is Gâteaux directional differentiable on $\mbox{bd}C$ satisfying that, for every point  $\ox\in \mbox{bd}C$, the following representations are satisfied\\
$$P^{'}_{C}(\ox)(w)=\begin{cases}
	w-\dfrac{1}{r^2}\la \ox-c,w\ra(\ox-c) \ \ \ \mbox{if} \ \ \ w\in\ox^\uparrow_{(c,r)};\\
	0 \ \ \ \ \ \ \ \ \ \ \ \ \ \ \ \ \ \ \ \ \ \ \ \ \ \ \ \ \ \ \ \mbox{if}\ \  \ w=\ox;\\w \ \ \ \  \ \ \ \ \ \ \ \ \ \ \ \ \ \ \ \ \ \ \ \ \ \  \ \ \ \  \mbox{if} \  \ \ \ w\in\ox^\downarrow_{(c,r)}.
\end{cases}$$  

(b)  $P_{C}$ is not Fréchet differentiable at any point $\ox\in \mbox{bd}C$. That is, $\nabla P_{C}(\ox)$ 
does not exist, for any $\ox\in \mbox{bd}C$. 

\end{Theorem}
{\bf Proof.} By Lemma \ref{Lem33}, we have $$P_C=P_{r\B}\circ\psi + h$$
$$P_{r\B}=P_C\circ\ph - h$$
where $\psi,\ph,h:H\to H$ with $\psi(x)=x-c, \ph(x)=x+c$, and $h(x)=c$ for all $x\in H.$ Note that $\psi,\ph,h$ are strictly differentiable at any $x\in H$ and $\nabla\psi(x)=\nabla\ph(x)=I_H, \ \nabla h(x)=0_H.$\\

(i). For any given $\ox\in C$, we have $\ox-c\in r\B$, by Proposition \ref{Pro32} and Propositions \ref{Pro23}, \ref{Pro24}, we get $P_C$ is strictly Fréchet differentiable at point $\ox$, and
$$\begin{array} {rl}\nabla P_C(\ox)&=\nabla P_{r\B}(\ox-c)\nabla\psi(\ox)+\nabla h(\ox)\\
&=I_H+0_H\\
&=I_H. \end{array}$$

(ii).  For $\ox\in H\backslash C$, we have $\ox-c\in H\backslash r\B$, by Proposition \ref{Pro32} and Propositions \ref{Pro23}, \ref{Pro24}, we get $P_C$ is strictly Fréchet differentiable at point $\ox$, and
$$\begin{array} {rl}\nabla P_C(\ox)(u)&=\nabla P_{r\B}(\ox-c)\big(\nabla\psi(\ox)(u)\big)+\nabla h(\ox)(u)\\
	&=\nabla P_{r\B}(\ox-c)(u)+\nabla h(\ox)(u)\\
	&=\dfrac{r}{\|\ox-c\|}\Big(u-\dfrac{\la \ox-c, u\ra}{\|\ox-c\|^2}(\ox-c)\Big)+0\\
	&=\dfrac{r}{\|\ox-c\|}\Big(u-\dfrac{\la \ox-c, u\ra}{\|\ox-c\|^2}(\ox-c)\Big),\ \ \mbox{for all} \ u\in H. \end{array}$$

(iii). Part (a) of (iii) follows from part (iii) in \cite[Theorem 5.2]{Li23}. Now we prove part (b) of (iii). For an arbitrary given $\ox\in \mbox{bd}C$, suppose to contrary that $P_C$ is Fréchet differentiable at $\ox$.  By Proposition \ref{Pro32} and Propositions \ref{Pro23}, \ref{Pro24}, $P_{r\B}=P_C\circ\ph-h$  is  Fréchet differentiable at $\ox-c\in\mbox{bd}r\B$. This contradicts the Proposition \ref{Pro32} (iii(b)).  So, $P_C$ is not Fréchet differentiable at $\ox.$                       

\hfill $\square$

To end this section, we give some examples to demonstrate the results of Theorem \ref{Thm34}.
\begin{Example} {\rm 
	Let $H=\R$, for any $c\in \R, r>0, \ C=\B(c,r)=[c-r;c+r],$ $C^o=(c-r;c+r)$ and $\mbox{bd}C=\{c-r;c+r\}$. For $x\in \R$, we have
	$$\nabla P_C(x)(u)=\begin{cases}
		s\ \ \ \ \ \ \ \ \ \ \ \ \ \ \ \ \ \ \ \ \ \mbox{if}\ \ x\in (c-r; c+r)\\
		0\ \ \ \ \ \ \ \ \ \ \ \ \ \ \ \ \ \ \  \ \ \mbox{if}\ \ x\not\in [c-r; c+r]\\
		\mbox{does not exits }\ \ \ \ \ \mbox{if}\ \ x=c-r \ \mbox{or}\ c+r
	\end{cases}\ \ \ \ \ \ \ \mbox{for every}\ u\in\R.$$
}
\end{Example}

\begin{Example} {\rm
		Let $H=\R^2$ and $C=\B(c,r)$ is the closed ball in $\R^2$ with radius $r>0$ and center $c=(c_1,c_2)\in\R^2$. For any $x=(x_1,x_2)\in \R^2,$ we have
	$$\begin{array} {rl} &\nabla P_C(x)(u)\\
		&=\begin{cases}
		u\ \ \ \ \ \ \ \ \ \ \ \ \ \ \ \ \ \ \ \ \ \ \ \ \ \ \ \ \ \ \ \ \ \ \ \ \ \ \ \ \ \ \ \ \ \ \ \ \ \ \ \  \ \ \ \ \ \ \ \ \mbox{if}\ \ (x_1-c_1)^2+(x_2-c_2)^2<r^2\\
		\dfrac{r}{\sqrt{(x_1-c_1)^2+(x_2-c_2)^2}}\bigg((u_1,u_2)-\dfrac{u_1(x_1-c_1)+u_2(x_2-c_2)}{(x_1-c_1)^2+(x_2-c_2)^2}(x_1-c_1,x_2-c_2) \bigg)\\  \ \ \ \ \ \ \ \ \ \ \ \ \ \ \ \ \ \ \ \ \  \ \ \ \ \ \ \ \ \ \ \ \ \ \ \ \ \ \ \ \ \ \ \ \ \ \ \ \ \ \ \ \ \ \ \ \ \ \ \  \ \ \mbox{if}\ \ (x_1-c_1)^2+(x_2-c_2)^2>r^2\\
		\mbox{does not exits }\ \ \ \ \ \ \ \ \ \ \ \ \ \ \ \ \ \ \ \ \ \ \ \ \ \ \ \ \ \ \ \ \ \ \ \ \ \ \ \ \  \ \ \ \mbox{if}\ \ (x_1-c_1)^2+(x_2-c_2)^2=r^2
	\end{cases}\\&  \mbox{for every}\ u=(u_1,u_2)\in\R.\end{array}$$
}
\end{Example}


\section{Strict Fréchet differentiability of the metric projection operator onto the second-order cone}\label{Sec3}
\setcounter{equation}{0}

In this section, we consider $\mathcal{K}$ is the m-dimensional second-order cone (also called Lorentz cone), i.e.,
$$\mathcal{K}:=\{x=(z_1,z_2)\in\R\times\R^{m-1}|\ z_1\geq\|z_2\|\}$$
with $m\geq 2.$ The topological interior and the boundary
of $\mathcal{K}$ are
$$\mathcal{K}^o=\{z=(z_1,z_2)\in\R\times\R^{m-1}|\ z_1>\|z_2\|\},\ \ \mbox{bd}\mathcal{K}=\{z=(z_1,z_2)\in\R\times\R^{m-1}|\ z_1=\|z_2\|\}.$$
For any given nonzero vector $z:=(z_1,z_2)\in\R\times\R^{m-1}$, we denote by
$$\lambda_i(z)=z_1+(-1)^i\|z_2\|\ \ \mbox{and} \ \ c_i(z)=\frac{1}{2}(1,(-1)^i\oz_2),\ \ \mbox{for} \  i=1,2,$$
where $$\oz_2=\begin{cases}
\dfrac{z_2}{\|z_2\|}\  \ \ \ \ \  \ \ \ \ \  \ \ \ \ \  \ \ \ \ \  \ \ \ \ \  \ \ \ \ \  \ \ \ \ \  \ \ \ \ \ \ \ \ \  \ \ \ \ \  \ \ \ \ \  \ \ \ \ \  \ \ \ \ \ \ \ \mbox{if}\ z_2\not=0,\\
w\ \ \ \ \ \ \mbox{with is any vector in} \ \R^{m-1}\ \mbox{satisfying}\ \|w\|=1,\ \ \ \mbox{if}\ z_2=0.
\end{cases}$$
Then we know from \cite{OS08} that $z$ has the following spectral decomposition
$$z=\lambda_1(z)c_1(z)+\lambda_2(z)c_2(z),$$
and a metric projecton operator on the second-order cone $P_\mathcal{K}$ can be written as
$$P_\mathcal{K}(z)=\lambda_1(z)_+c_1(z)+\lambda_2(z)_+c_2(z).$$

The following result is given in \cite{OS08, YZ18}, which proves the Fréchet differentiability of $P_\mathcal{K}$.
\begin{Proposition} \label{Pro41} 
	Let $\mathcal{K}$ be the m-dimensional second-order cone. Then, the metric projection operator  $P_{\mathcal{K}}:\R^m\to\mathcal{K}$ has the following differentiability properties.
	
	(i) $P_\mathcal{K} $  is  Fréchet differentiable on $\mathcal{K}^o$, and $$\nabla P_{\mathcal{K}}(z)=I_{\R^m}\ \ \mbox{for every}\ z\in \mathcal{K}^o.$$
	
	(ii)  $P_\mathcal{K} $  is  Fréchet differentiable on  $ -\mathcal{K}^o$, and $$\nabla P_{\mathcal{K}}(z)=0_{\R^m} \ \ \mbox{for every}\ z\in -\mathcal{K}^o,$$
	where $0_{\R^m}:\R^m\to\R^m$ with $0_{\R^m}(u)=0$ for all $u\in \R^m.$
	
	(iii) $P_\mathcal{K} $  is  Fréchet differentiable on  $\R^m\backslash(-\mathcal{K}\cup\mathcal{K})$, and for every $z\in \R^m\backslash(-\mathcal{K}\cup\mathcal{K})$,   $$\nabla P_\mathcal{K}(z)=\dfrac{1}{2}\Big(1+\dfrac{z_1}{\|z_2\|}\Big)I_{\R^m}+\dfrac{1}{2}\begin{bmatrix}
	-\dfrac{z_1}{\|z_2\|}&\oz_2^T\\
	\oz_2&-\dfrac{z_1}{\|z_2\|}\oz_2\oz_2^T
	\end{bmatrix}.$$

(iv) $P_\mathcal{K} $  is Gâteaux directionally differentiable at any $z\in \mbox{bd}\mathcal{K}$ and for any $h\in\R^m,$
$$P^{'}_\mathcal{K}(z)(h)=\begin{cases}
h-2(c_1(z)^Th)_-c_1(z)\ \ \ \ \ \ \ \ \ \ \ \ \ \ \mbox{if}\ \ z\in \mbox{bd}\mathcal{K}\backslash\{0\},\\
2(c_2(z)^Th)_+c_2(z)\ \ \ \ \ \ \ \ \ \ \ \ \ \ \ \ \ \ \ \mbox{if}\ \ z\in \mbox{bd}\mathcal{K}\backslash\{0\},\\
(\lambda_1(h))_+c_1(h)+(\lambda_2(h))_+c_2(h)\ \ \ \mbox{if}\ \  z=0.
\end{cases}$$
\end{Proposition}

In the following theorem, we prove the  strict Fréchet differentiability of the metric projection operator $P_\mathcal{K}$ on $\R^m\backslash(-\mbox{bd}\mathcal{K}\cup\mbox{bd}\mathcal{K})$ and show that $P_\mathcal{K}$ is not differentiable at $z\in -\mbox{bd}\mathcal{K}\cup\mbox{bd}\mathcal{K}$.
\begin{Theorem} \label{Thm42}

	Let $\mathcal{K}$ be the m-dimensional second-order cone. Then, the metric projection operator $P_{\mathcal{K}}:\R^m\to\mathcal{K}$ has the following differentiability properties.
	
	(i) $P_\mathcal{K} $  is strictly Fréchet differentiable on $\mathcal{K}^o$, and $$\nabla P_{\mathcal{K}}(z)=I_{\R^m} \ \ \mbox{for every}\ z\in \mathcal{K}^o.$$
	
	(ii)  $P_\mathcal{K} $  is strictly Fréchet differentiable on  $ -\mathcal{K}^o$, and $$\nabla P_{\mathcal{K}}(z)=0_{\R^m}\ \ \mbox{for every}\ z\in -\mathcal{K}^o.$$
	
	(iii) $P_\mathcal{K} $  is strictly Fréchet differentiable on  $\R^m\backslash(-\mathcal{K}\cup\mathcal{K})$, and for any $z\in \R^m\backslash(-\mathcal{K}\cup\mathcal{K})$,   $$\nabla P_\mathcal{K}(z)=\dfrac{1}{2}\Big(1+\dfrac{z_1}{\|z_2\|}\Big)I_{\R^m}+\dfrac{1}{2}\begin{bmatrix}
		-\dfrac{z_1}{\|z_2\|}&\oz_2^T\\
		\oz_2&-\dfrac{z_1}{\|z_2\|}\oz_2\oz_2^T
	\end{bmatrix}.$$
	
	(iv) For the subset $-\mbox{bd}\mathcal{K}\cup\mbox{bd}\mathcal{K}$, we have
	
\ \ \ 	(a) $P_\mathcal{K} $  is Gâteaux directionally differentiable at any $z\in-\mbox{bd}\mathcal{K}\cup \mbox{bd}\mathcal{K}$ and for any $h\in\R^m,$
	$$P^{'}_\mathcal{K}(z)(h)=\begin{cases}
		h-2(c_1(z)^Th)_-c_1(z)\ \ \ \ \ \ \ \ \ \ \ \ \ \ \mbox{if}\ \ z\in \mbox{bd}\mathcal{K}\backslash\{0\},\\
		2(c_2(z)^Th)_+c_2(z)\ \ \ \ \ \ \ \ \ \ \ \ \ \ \ \ \ \ \ \mbox{if}\ \ z\in -\mbox{bd}\mathcal{K}\backslash\{0\},\\
		(\lambda_1(h))_+c_1(h)+(\lambda_2(h))_+c_2(h)\ \ \ \mbox{if}\ \  z=0.
	\end{cases}$$

\ \ \ (b) $P_\mathcal{K} $ is not Fréchet differentiable at any point $z\in-\mbox{bd}\mathcal{K}\cup \mbox{bd}\mathcal{K}$.
\end{Theorem}
{\bf Proof.}
(i) For any $z\in\mathcal{K}^o$,  there is a positive number $\delta$ such that 
$$\B(z,\delta)\subset\mathcal{K}^o.$$
We have $P_\mathcal{K}(u)=u,\ \ \mbox{for every}\ u\in \B(z,\delta).$\\
This implies that
$$\begin{array}{rl}
	\lim\limits_{(u,v)\to(z,z)}\dfrac{P_\mathcal{K}(u)-P_\mathcal{K}(v)-I_{\R^m}(u-v)}{\|u-v\|}=\lim\limits_{(u,v)\to(z,z)}\dfrac{u-v-(u-v)}{\|u-v\|}=0.
\end{array}$$
Hence, $P_\mathcal{K} $  is strictly Fréchet differentiable at $z\in \mathcal{K}^o$ with $\nabla P_{\mathcal{K}}(z)=I_{\R^m}$.\\

(ii)  Similar to proof of part (i), for any $z\in-\mathcal{K}^o$,  there exists $\delta>0$ such that $\B(z,\delta)\subset-\mathcal{K}^o$, and $P_\mathcal{K}(u)=0,\ \ \mbox{for every}\ u\in \B(z,\delta).$\\
This implies that
$$\begin{array}{rl}
	\lim\limits_{(u,v)\to(z,z)}\dfrac{P_\mathcal{K}(u)-P_\mathcal{K}(v)-0_{\R^m}(u-v)}{\|u-v\|}=\lim\limits_{(u,v)\to(z,z)}\dfrac{0-0-0}{\|u-v\|}=0.
\end{array}$$
So, $P_\mathcal{K} $  is strictly Fréchet differentiable at $z\in \mathcal{K}^o$ with $\nabla P_{\mathcal{K}}(z)=0_{\R^m}$.\\

(iii) For any $z\in \R^m\backslash(-\mathcal{K}\cup\mathcal{K})$, by Proposition \ref{Pro41}, there exists $\delta=\delta_1\|z\|>0$ such that $\B(z,\delta)\subset \R^m\backslash(-\mathcal{K}\cup\mathcal{K})$  and $P_\mathcal{K} $  is Fréchet differentiable at every $u\in \B(z,\delta)$ with $$\nabla P_\mathcal{K}(u)=\dfrac{1}{2}\Big(1+\dfrac{u_1}{\|u_2\|}\Big)I_{\R^m}+\dfrac{1}{2}\begin{bmatrix}
	-\dfrac{u_1}{\|u_2\|}&\ou_2^T\\
	\ou_2&-\dfrac{u_1}{\|u_2\|}\ou_2\ou_2^T
\end{bmatrix}.$$
For any $u,v\in \B(z,\delta)$, we have
$$\begin{array} {rl}
	\nabla P_\mathcal{K}(u)-\nabla P_\mathcal{K}(v)&=
	\dfrac{1}{2}\Bigg(\dfrac{u_1}{\|u_2\|}-\dfrac{v_1}{\|v_2\|}\Bigg)I_{\R^m}+\dfrac{1}{2}\begin{bmatrix}
		\dfrac{v_1}{\|v_2\|}-\dfrac{u_1}{\|u_2\|}&\ou_2^T-\ov_2^T\\
		\ou_2-\ov_2&\dfrac{v_1}{\|v_2\|}\ov_2\ov_2^T-\dfrac{u_1}{\|u_2\|}\ou_2\ou_2^T
	\end{bmatrix}\\
&=\dfrac{1}{2}\dfrac{u_1\|v_2\|-v_1\|u_2\|}{\|u_2\|\|v_2\|}I_{\R^m}
+\dfrac{1}{2}\begin{bmatrix}
	\dfrac{v_1\|u_2\|-u_1\|v_2\|}{\|v_2\|\|u_2\|}&\ou_2^T-\ov_2^T\\
	\ou_2-\ov_2&\dfrac{v_1}{\|v_2\|}\ov_2\ov_2^T-\dfrac{u_1}{\|u_2\|}\ou_2\ou_2^T
\end{bmatrix}\\
\end{array}$$
   Consider matrix $$A=\begin{bmatrix}
	\dfrac{v_1\|u_2\|-u_1\|v_2\|}{\|v_2\|\|u_2\|}&\ou_2^T-\ov_2^T\\
	\ou_2-\ov_2&\dfrac{v_1}{\|v_2\|}\ov_2\ov_2^T-\dfrac{u_1}{\|u_2\|}\ou_2\ou_2^T
\end{bmatrix}=[a_{ij}]_{m\times m}.$$
We have $$\begin{array}{rl} |a_{11}|&=\dfrac{\big|\|u_2\|(v_1-u_1)+u_1(\|u_2\|-\|v_2\|)\big|}{\|v_2\|\|u_2\|}
\leq\dfrac{\|u_2\|\big|v_1-u_1\big|+|u_1|\big|\|u_2\|-\|v_2\|\big|}{\|v_2\|\|u_2\|}\\
&\leq\dfrac{\|v_1-u_1\|+\|u_2-v_2\|}{\|v_2\|}\leq \dfrac{2\|u-v\|}{\|v_2\|}
\leq \dfrac{2\sqrt{2}}{\|v\|}\|u-v\|\leq \dfrac{2\sqrt{2}}{\|z\|(1-\delta_1)}\|u-v\|.\end{array}$$
$$\begin{array}{rl}|a_{j1}|&=|a_{1j}|=\Bigg|\dfrac{u_{2j}}{\|u_{2}\|}-\dfrac{v_{2j}}{\|v_2\|}\Bigg| =\dfrac{\big|u_{2j}\|v_2\|-v_{2j}\|u_2\|\big|}{|v_{2}\|\|u_{2}\|}\\&
\leq\dfrac{|u_{2j}|\big|\|v_2-\|u_2\|\big|+|u_{2j}-v_{2j}|\|u_2\|}{\|v_2\|\|u_2\|}		\leq\dfrac{\big|\|v_2-\|u_2\|\big|+|u_{2j}-v_{2j}|}{\|v_{2}\|}
	\\&		\leq \dfrac{2\|u-v\|}{\|v_2\|}
	\leq \dfrac{2\sqrt{2}}{\|z\|(1-\delta_1)}\|u-v\|
	,\ \ \mbox{for all}\ j=2,...,m.\end{array}$$

$$\begin{array}{rl}|a_{ij}|&=\Bigg|\dfrac{v_1v_{2i}v_{2j}}{\|v_2\|^3}-\dfrac{u_1u_{2i}u_{2j}}{\|u_2\|^3}\Bigg|\\
	&=\dfrac{\big|v_1v_{2i}v_{2j}\|u_2\|^3-u_1u_{2i}u_{2j}\|v_2\|^3\big|}{\|v_2\|^3\|u_2\|^3}\\
	&\leq \dfrac{\big| v_1v_{2i}v_{2j}-u_1u_{2i}u_{2j}\big|\|u_2\|^3+|u_1u_{2i}u_{2j}|\big|\|u_2\|^3-\|v_2\|^3\big|}{\|v_2\|^3\|u_2\|^3}\\
	&\leq \dfrac{\big| v_1v_{2i}v_{2j}-u_1u_{2i}u_{2j}\big|+\big|\|u_2\|^3-\|v_2\|^3\big|}{\|v_2\|^3}\\
	&=\dfrac{\big| v_1v_{2i}(v_{2j}-u_{2j})+v_1u_{2j}(v_{2i}-u_{2i})+u_{2i}u_{2j}(v_1-u_1)\big|+\big|\|u_2\|^3-\|v_2\|^3\big|}{\|v_2\|^3}\\
	&\leq\dfrac{\|v\|^2\|v-u\|+\|v\|\|u\|\|v-u\|+\|u\|^2\|v-u\|+\big|\|u_2\|^3-\|v_2\|^3\big|}{\frac{1}{2\sqrt{2}}\|v\|^3}\\
	&\leq\dfrac{2\sqrt{2}(\|v\|^2+\|v\|\|u\|+\|u\|^2)}{\|v\|^3}\|u-v\|+\dfrac{2\sqrt{2}\|u_2-v_2\|(\|u_2\|^2+\|u_2\|\|v_2\|+\|v_2\|^2)}{\|v\|^3}\\
	&\leq\dfrac{4\sqrt{2}(\|v\|^2+\|v\|\|u\|+\|u\|^2)}{\|v\|^3}\|u-v\|
	\\
	&\leq\dfrac{12\sqrt{2}(1+\delta_1)^2}{\|z\|(1-\delta_1)^3}\|u-v\|, \ \ \mbox{for all}\ \ i,j\in\{2,...,m\}.	
\end{array}$$
Where $u=(u_1,u_2), v=(v_1,v_2)\in \R^m$ with $u_2:=(u_{22},...,u_{2m}) ,v_2:=(v_{22},...,v_{2m}) \in\R^{m-1}$, notice that for any $u, v\in \B(z,\delta)$, we have $|u_1|<\|u_2\|, |v_1|<\|v_2\|$ and $(1-\delta_1)\|z\|\leq \|u\|, \|v\|\leq \|z\|(1+\delta_1)$.\\
This implies that $$\|A\|\leq\sqrt{\sum_{i,j=1}^{m}a_{ij}^2}\leq \dfrac{12\sqrt{2}m(1+\delta_1)^2}{\|z\|(1-\delta_1)^3}\|u-v\|. $$
Hence  
$$\begin{array}{rl}\|\nabla P_\mathcal{K}(u)-\nabla P_\mathcal{K}(v)\|&\leq \dfrac{1}{2}\Bigg|\dfrac{u_1\|v_2\|-v_1\|u_2\|}{\|u_2\|\|v_2\|}\Bigg|+\dfrac{1}{2}\|A\|\\
&\leq \dfrac{1}{2}\dfrac{\Big||u_1|\big|\|v_2\|-\|u_2\|\big|+\|u_2\||u_1-v_1|\Big|}{\|u_2\|\|v_2\|}+\dfrac{1}{2}\|A\|\\	
&\leq \dfrac{1}{2}\dfrac{\Big|\big|\|v_2\|-\|u_2\|\big|+|u_1-v_1|\Big|}{\|v_2\|}+\dfrac{1}{2}\|A\|\\		
	
&\leq \dfrac{1}{2}\dfrac{2\|v-u\|}{\frac{1}{\sqrt{2}}\|v\|}+\dfrac{1}{2}\dfrac{12\sqrt{2}m(1+\delta_1)^2}{\|z\|(1-\delta_1)^3}\|u-v\|\\
&\leq \dfrac{\sqrt{2}\|u-v\|}{\|z\|(1-\delta_1)}+\dfrac{6\sqrt{2}m(1+\delta_1)^2}{\|z\|(1-\delta_1)^3}\|u-v\|\\&
=\dfrac{\sqrt{2}(1-\delta_1)^2+6\sqrt{2}m(1+\delta_1)^2}{\|z\|(1-\delta_1)^3}\|u-v\|.
\end{array}$$
Therefore $\nabla P_\mathcal{K}(.)$ is Lipschitz continuous on $\B(z, \delta)$. By Lemma \ref{Pro31}, $P_{\mathcal{K}}$ is strict Fréchet differentiable at $z$.

(iv) Part (a) of (iv) follows from part (iv) in Proposition \ref{Pro41}. To prove part (b), we consider the following cases.

{\it Case 1.} $z=(z_1,z_2)\in\mbox{bd}\mathcal{K}\backslash\{0\}.$ Suppose to contrary that $P_\mathcal{K}$ is Fréchet differentiable at $z$. Then, there is a linear continuous mapping $\nabla P_\mathcal{K}(z):\R^m\to\R^m$, such that
$$\nabla P_\mathcal{K}(z)(u)=P^{'}(z)(u)\ \ \mbox{for all}\ \ u\in\R^m.$$
Choosing  $\widetilde z=(-z_1,z_2)$, by (iv.a), we have 
\begin{equation}\label{312}\nabla P_\mathcal{K}(z)(\widetilde z)=P^{'}(z)(\widetilde z)=\widetilde z-2\big(c_1(z)^T\widetilde z\big)_-c_1(z)=\widetilde z-\widetilde z=0.
\end{equation}
And 
$$\nabla P_\mathcal{K}(z)(-\widetilde z)=P^{'}(z)(-\widetilde z)=-\widetilde z-2\big(c_1(z)^T(-\widetilde z)\big)_-c_1(z)=-\widetilde z-0=-\widetilde z.$$
Thanks for the linearity of $\nabla P_\mathcal{K}$, we get $\nabla P_\mathcal{K}(z)(\widetilde z)=\widetilde z.$
This contradicts to \eqref{312}, which proves that $P_\mathcal{K}$ is not Fréchet differentiable at any point $z\in\mbox{bd}\mathcal{K}\backslash\{0\}.$
 
 {\it Case 2.} $z=(z_1,z_2)\in-\mbox{bd}\mathcal{K}\backslash\{0\}$, that is $-z_1=\|z_2\|\not=0.$ Suppose to contrary that $P_\mathcal{K}$ is Fréchet differentiable at $z$. Then, 
 $$\nabla P_\mathcal{K}(z)(u)=P^{'}(z)(u)\ \ \mbox{for all}\ \ u\in\R^m.$$
 Choosing  $\widetilde z=(-z_1,z_2)$, by (iv.a), we have 
 \begin{equation}\label{111}\nabla P_\mathcal{K}(z)(\widetilde z)=P^{'}(z)(\widetilde z)=2\big(c_2(z)^T\widetilde z\big)_+c_2(z)=\widetilde z.
 \end{equation}
 And 
 $$\nabla P_\mathcal{K}(z)(-\widetilde z)=P^{'}(z)(-\widetilde z)=2\big(c_2(z)^T(-\widetilde z)\big)_+c_2(z)=0.$$
 Thanks for the linearity of $\nabla P_\mathcal{K}$, we get $\nabla P_\mathcal{K}(z)(\widetilde z)=0.$
 This contradicts to \eqref{111}, which proves that $P_\mathcal{K}$ is not Fréchet differentiable at any point $z\in-\mbox{bd}\mathcal{K}\backslash\{0\}.$ 

{\it Case 3.} For $z=0\in\R^m$. Suppose to contrary that $P_\mathcal{K}$ is Fréchet differentiable at $z$. Then, choosing $\widehat z=(1,0)\in\R\times\R^{m-1},$  by (iv.a), we have 
\begin{equation}\label{222}\nabla P_\mathcal{K}(z)(\widehat z)=P^{'}(z)(\widehat z)=\big(\lambda_1(\widehat z)\big)_+c_1(\widehat z)+\big(\lambda_2(\widehat z)\big)_+c_2(\widehat z)=c_1(\widehat z)+c_2(\widehat z)=\widehat z.
\end{equation}
And 
$$\nabla P_\mathcal{K}(z)(-\widehat z)=P^{'}(z)(-\widehat z)=\big(\lambda_1(-\widehat z)\big)_+c_1(-\widehat z)+\big(\lambda_2(-\widehat z)\big)_+c_2(-\widehat z)=0.$$
Thanks for the linearity of $\nabla P_\mathcal{K}$, we have $\nabla P_\mathcal{K}(z)(\widehat z)=0.$ This contradicts to \eqref{222}, which proves that $P_\mathcal{K}$ is not Fréchet differentiable at $z=0$.
$\hfill\square$

\begin{Example}{\rm 
Let $\mathcal{K}=\{z=(z_1,z_2,z_3)\in\R^3|z_1\geq\|z_2\|\}$ is a second-order in $\R^3$. Then, we have

(i) If $z_1>\sqrt{z_2^2+z_3^2}$, then 
$$\nabla P_\mathcal{K}(z)(h)=h, \ \ \ \mbox{for any}\ h\in \R^3.$$

(ii) If $-z_1>\sqrt{z_2^2+z_3^2}$, then 
$$\nabla P_\mathcal{K}(z)(h)=0, \ \ \ \mbox{for any}\ h\in \R^3.$$

(iii) If $|z_1|<\sqrt{z_2^2+z_3^2}$, then 
$$\nabla P_\mathcal{K}(z)(h)=\dfrac{1}{2\sqrt{(z_2^2+z_3^2)^3}}\begin{bmatrix}
	\sqrt{(z_2^2+z_3^2)^3}&z_2(z_2^2+z_3^2)&z_3(z_2^2+z_3^2)\\
	z_2(z_2^2+z_3^2)&	\sqrt{(z_2^2+z_3^2)^3}+z_1z_3^2&-z_1z_2z_3\\
	z_3(z_2^2+z_3^2)&-z_1z_2z_3& \sqrt{(z_2^2+z_3^2)^3}+z_1^2z_2^2
\end{bmatrix}\begin{bmatrix}
h_1\\h_2\\h_3
\end{bmatrix} $$$\ \ \ \mbox{for any}\ h=(h_1, h_2, h_3)\in \R^3.$

(iv)  If $|z_1|=\sqrt{z_2^2+z_3^2}$, one has $\nabla P_\mathcal{K}(z)$ does not exist.
}
\end{Example}

\section{Concluding remarks}

The main result of this paper is to provide a simpler proof for the result in \cite{Li24}, which mainly proving strict Fréchet differentiability of the metric projection operator onto closed balls centered at origin, and we generalize this result for the metric projection operator onto closed balls with center at arbitrarily given point $c$ in Hilbert spaces. In addition, we also study the strict differentiability of the metric projection operator onto the second-order cones in Euclidean spaces. In the near future, we are going to conduct more researches on the strict Fréchet differentiability of $P_C$, with $C$ being special a subset in Hilbert spaces.

.
\section*{Acknowledgments}
The author would like to thank Vietnam Institute for Advanced Study in Mathematics for hospitality during her post-doctoral fellowship of the Institute in 2022--2023.

\end{document}